\documentclass[twoside,11pt]{article}
\usepackage{graphicx,amscd,amsfonts,amsmath,amsthm,amssymb,latexsym,amsfonts,color}
\usepackage{epsfig,float,epstopdf,array,bm,cite}
\usepackage{cases}
\usepackage{array}
\usepackage{multirow}
\usepackage[table]{xcolor}
\usepackage{relsize}
\usepackage{graphicx}
\usepackage{float}
\usepackage[bookmarksnumbered, colorlinks,plainpages]{hyperref}

\def\diag{\textrm{diag}}

\def\h{\hspace{-0.2cm}}
 \newcommand{\cred}[1]{{\color{black}  #1}}
  \newcommand{\cblu}[1]{{\color{black}  #1}}

\newtheorem{theorem}{\bf Theorem}
\newtheorem{lemma}{\bf Lemma}

\newtheorem{remark}{\bf Remark}
\def\G{\mathcal{G}}
\def\K{\mathcal{K}}
\def\M{\mathcal{M}}
\def\T{\mathcal{T}}
\def\P{\mathcal{P}}
\def\Q{\mathcal{Q}}
\def\B{\mathcal{B}}
\def\R{\mathcal{R}}
\def\S{\mathcal{S}}
\def\h{\hspace{-0.2cm}}

\begin{document}
\title{\bf A new iterative method for solving a class of  two-by-two block complex linear systems}
\author{\small\bf Davod Khojasteh Salkuyeh\thanks{\noindent Corresponding author. \newline
		Emails:	khojasteh@guilan.ac.ir} \\[2mm]
	\textit{{\small Faculty of Mathematical Sciences, University of Guilan, Rasht, Iran}} \\
	\textit{{\small Center of Excellence for Mathematical Modelling, Optimization and Combinational}}\\
	\textit{{\small Computing (MMOCC), University of Guilan, Rasht, Iran}}
}
\date{}
\maketitle
\vspace{-0.5cm}

\noindent\hrulefill\\
{\bf Abstract.}  We present a stationary iteration method, namely \cblu{Alternating Symmetric positive definite and  Scaled symmetric positive semidefinite Splitting (ASSS)}, for solving the system of linear equations obtained by using finite element discretization of a distributed optimal control problem together with  time-periodic parabolic equations. An upper bound for the spectral radius of the iteration method is given which is always 
\cblu{less than 1}. So  convergence of the ASSS iteration method is guaranteed. The induced ASSS preconditioner is applied to accelerate the convergence speed of the GMRES method for solving the system. Numerical results are presented to demonstrate the effectiveness of both the ASSS iteration method and the ASSS preconditioner.    
 \\[-3mm]


\noindent\hrulefill

\pagestyle{myheadings}\markboth{D.K. Salkuyeh  }{A new iterative method for solving a class of  two-by-two block complex linear systems}
\thispagestyle{empty}
\section{Introduction}\label{Sec1}

Consider the distributed  control problems of the form: (see \cite{Lions,Kollmann,ZhengCAM,ZhengAML})
\begin{eqnarray*}
	\min_{y,u}~~~ \frac{1}{2}\int_0^T\int_\Omega  |y(x,t)-y_d(x,t)|^2dxdt &\h+\h& \frac{\nu}{2}\int_0^T\int_\Omega |u(x,t)|^2dxdt,\label{Eq001}\\
	\hspace{2cm}\text{s.t.}~~~~~~~~~~\frac{\partial}{\partial t}y(x,t)-\Delta y(x,t)&\h=\h&u(x,t)~~ {\rm in}~~ Q_T,~~~~~~~~~~~~~~~~~~~\hspace{2cm}\label{Eq002}\\
	y(x,t)&\h=\h&0 ~~{\rm on}~~ \Sigma_T,~~~~~~~~~~~~~~~~~~~~~~~~~~~~~~~~~~\hspace{0.1cm}\hspace{2cm}\label{Eq0003}\\
	y(x,0)&\h=\h&y(x,T) ~~{\rm on} ~~\partial \Omega,~~~~~~~~~~~~~~~~~~~~~~~~~~~\hspace{2cm}\label{Eq0004}\\
	u(x,0)&\h=\h&u(x,T)~~ {\rm in}~~ \Omega ,~~~~~~~~~~~~~~~~~~~~~~~~~~~~~\hspace{0.05cm}\hspace{2cm}\label{Eq0005}
\end{eqnarray*}
where $\Omega$ is  an open and bounded domain in $\mathbb{R}^d$ ($d\in\{1,2,3\})$ and its boundary $\partial \Omega $  is Lipschitz-continuous. We introduce the space-time cylinder $Q_T=\Omega \times (0,T)$ and its lateral surface \cblu{$\Sigma_T=\partial\Omega\times (0,T)$}. Here, $\nu$ is a regularization parameter,  $y_d(x,t)$ is a desired state and
$T>0$. 
Problems of the form above arise widely in many areas of science and engineering \cite{Kollmann,Krendl}.

We may assume that $y_d(x,t)$ is time-harmonic, i.e.,
$y_d(x,t)=y_d(x)e^{i\omega t}$, with $\omega={2\pi k}/{T}$ for some $k\in\Bbb{Z}$.
If we substitute $y(x,t)$, $y_d(x,t)$ and $u(x,t)$ in the problem then we get the following time-independent problem   
\begin{eqnarray*}
	\min_{y,u}  && \frac{1}{2}\int_\Omega  |y(x)-y_d(x)|^2dx +\frac{\nu}{2}\int_\Omega |u(x)|^2dx,\\
	s.t:        && \cblu{i\omega y(x)}-\Delta y(x)=u(x)~~\text{in}~~ \Omega,~~~~~~~~~~~~~~~~~~~\label{EqMain}\\
	&& \hspace{1.9cm} \cblu{y(x)=0}, \hspace{0.6cm} \text{on}~~ \partial\Omega.
\end{eqnarray*}
We assume that $V_h\subset H_{0}^1(\Omega)$ is an $m$-dimensional vector space spanned by the basis\linebreak $\{\varphi_1,\varphi_2,\ldots,\varphi_m\}$. The subspace  $V_h$ is utilized for computing both the functions $y$ and $u$.
By applying the approach of discretize-then-optimization (see \cite{Rees}), the above problem can be described as the following form  \cite{Krendl}   
\begin{equation}\label{Eqqa}
\begin{array}{lll}
	\displaystyle\min_{y,u}  && \frac{1}{2}  (\bar{y}-\bar{y}_d)^*M (\bar{y}-\bar{y}_d) +\frac{\nu}{2} \bar{u}^*M \bar{u},\\
	s.t:        && \cblu{iwM} \bar{y}+K \bar{y}=M\bar{u},
\end{array}
\end{equation}
where the discretized negative Laplacian  is interpreted by the matrix $K=(\int_{\Omega} \nabla\varphi_i . \nabla\varphi_j dx) \in\Bbb{R}^{m\times m}$ (the stiffness matrix) and  $M=(\int_{\Omega} \varphi_i \varphi_j dx) \in\Bbb{R}^{m\times m}$ denotes the mass matrix. Here, the vectors $\bar{y}$, $\bar{y}_d$, and $\bar{u}$ denote the coefficients  of the basis  functions in $V_h$. We define the Lagrangian functional for the discretized problem as
\[
{\cal L}(\bar{y},\bar{u},\bar{p})=\frac{1}{2}  (\bar{y}-\bar{y}_d)^*M (\bar{y}-\bar{y}_d) +\frac{\nu}{2} \bar{u}^*M \bar{u}+
{\bar p}^*(\cblu{i\omega M} \bar{y}+K \bar{y}-M\bar{u}),
\]
with $\bar{p}$ being the Lagrange multiplier of the constraint. By using the  Lagrange multiplier technique, we set $\nabla {\cal L}(\bar{y},\bar{u},\bar{p})=0$, which is equivalent to
\begin{equation}\label{MainSystem}
\begin{pmatrix}
M & 0 & K-i\omega M \\
0        & \nu M & -M \\
K+i\omega M    &-M & 0
\end{pmatrix}
\begin{pmatrix}
{\bar{y}} \\
{ \bar{u}}\\
{\bar{p}}
\end{pmatrix}
=
\begin{pmatrix}
M{\bar{y}_d} \\
0\\
0
\end{pmatrix}.
\end{equation}
It is known that, Eq. \eqref{MainSystem} gives the necessary and sufficient  conditions for the existence of a solution  for the problem \eqref{Eqqa}. It follows from the second equation that $\bar{u}=\frac{1}{\nu}\bar{p}$ and substituting $\bar{u}$ in the third equation, gives the following complex system of equations 
\begin{equation*}\label{Eq008}
\left\{
\begin{array}{rl}
 M\bar{y}+(K-i{\omega}M)\bar{p}&=M{\bar{y}_d},\\
 (K+i{\omega}M)\bar{y}-\frac{1}{\nu}M\bar{p}&=0. 
\end{array} \right.
\end{equation*}
The above system can be equivalently rewritten as  

\begin{equation}\label{MainSystem15}
A{x} =
\begin{pmatrix}
M & \sqrt{\nu}(K-i\omega M) \\
\sqrt{\nu}(K+i\omega M)      & -M
\end{pmatrix}
\begin{pmatrix}
{ \bar{y}} \\
{ \bar{q}}
\end{pmatrix}
=
\begin{pmatrix}
{\hat{y}_d} \\
0
\end{pmatrix}
=b,
\end{equation}
where $\bar{q}=\bar{p}/\sqrt{\nu}$ and $\hat{y}_d=M{\bar{y}_d}$.

Since the systems of the form \eqref{MainSystem15} are  of very large size, it is important to employ iterative methods in incorporated with suitable preconditioners for solving these systems. 
In \cite{Krendl}, Krendl et al.  presented the real block diagonal ($\mathcal{P}_{BD}$) and the alternative indefinite ($\mathcal{P}_{AI}$) preconditioners 
\begin{eqnarray}
&&\mathcal{P}_{BD}=
\begin{pmatrix}
 M+\sqrt{\nu}(K+\omega M) & 0 \\
0      & M+\sqrt{\nu}(K+\omega M)
\end{pmatrix},\label{Eqq03}\\
&&\mathcal{P}_{AI}=
\begin{pmatrix}
0 & M+\sqrt{\nu}(K-i\omega M)\\
M+\sqrt{\nu}(K+i\omega M)    & 0
\end{pmatrix},\label{Eqq04}
\end{eqnarray}
for the system \eqref{MainSystem15}. Zheng et al. in \cite{ZhengAML} presented an improved version of the block-diagonal preconditioner $\mathcal{P}_{BD}$. Zheng et al. in \cite{ZhengCAM} designed  the block alternating splitting (BAS) iteration method for solving the system \eqref{MainSystem15} which can be written  as
\begin{equation*}\label{Eq035}
\left\{
\begin{array}{rl}
(\alpha V+H_1)x^{(k+\frac{1}{2})}&=(\alpha V-S_1)x^{(k)}+\mathcal{P}_1 b,~~~~~~ \\
(\alpha V+H_2)x^{(k+1)}&=(\alpha V-S_2)x^{(k+\frac{1}{2})}+\mathcal{P}_2 b, ~~~
\end{array} \right.
\end{equation*}
where $\alpha>0$, $V$ is a symmetric positive definite matrix,
\[
H_1=\begin{pmatrix}
M& 0 \\
0& M
\end{pmatrix},\quad
H_2=\begin{pmatrix}
\sqrt{\nu}K& 0 \\
0& \sqrt{\nu}K
\end{pmatrix}
\]
and
\cblu{\begin{eqnarray*}
S_1 &=& \frac{1}{1+\omega^2 \nu}
\begin{pmatrix}
-{i\omega \nu }K &  {\sqrt{\nu}}K \\
-{\sqrt{\nu}}K   & {i\omega \nu }K
\end{pmatrix},\qquad
S_2 =\begin{pmatrix}
i\sqrt{\nu}\omega M &  -M \\
M & -i\sqrt{\nu}\omega M
\end{pmatrix},\\
\mathcal{P}_1&=&\frac{1}{1+\omega^2 \nu}
\begin{pmatrix}
-I &  -i\omega \sqrt{\nu} I \\
i \omega \sqrt{\nu} I   & -I
\end{pmatrix},\qquad
\mathcal{P}_2=\frac{1}{1+\omega^2 \nu}
\begin{pmatrix}
0 & I \\
I & 0
\end{pmatrix}.
\end{eqnarray*}}
\cred{They proposed using the matrix $V=H_1$ as the preconditioner in practical implementation}. They proved that the BAS iteration method is convergent if $\alpha\geq \nu \omega^2/2$. Numerical results presented in \cite{ZhengCAM} demonstrate that the BAS iteration method outperforms the GMRES method \cite{SaadBook,GMRES}. On the other hand, the BAS iteration method induces the preconditioner
\[
P_{BAS}(\alpha)=(1+\alpha)P(\alpha) \begin{pmatrix}
	\alpha M+\sqrt{\nu} K  & 0 \\
	                    0  & \alpha M+\sqrt{\nu} K
	\end{pmatrix},
\] 
with 
\cblu{\[
 P(\alpha)=\frac{1}{\alpha(2+\omega^2\nu)} \begin{pmatrix}
I &  (1+\omega^2\nu-i\omega\sqrt{\nu})I \\
(1+\omega^2\nu+i\omega\sqrt{\nu})I  & -I
\end{pmatrix}.
\]}
They also numerically illustrated that the parameter $\alpha=1+\nu \omega^2$ usually gives good results for the BAS iteration method and 
$\alpha=(1+\nu \omega^2)/(1+\sqrt{\nu}\omega)$ is a good choice for the BAS preconditioner. \cred{The reported numerical results in \cite{ZhengCAM} showed that the BAS preconditioner is superior to the preconditioners $\mathcal{P}_{BD}$ and $\mathcal{P}_{AI}$.}  

\cred{Recently, Mirchi and Salkuyeh \cite{SalkuyehMATCOM} stated the single block splitting (SBS) iteration method for solving  \eqref{MainSystem15}. They proved that if 
\begin{equation*}\label{25006}
\frac{\nu \omega\mu_{\max}}{\sqrt{(1+\omega^2\nu )^2+\nu \mu_{\max}^{2}}}<1,
\end{equation*}
then the SBS method is convergent, where $\mu_{\max} $ denotes the largest eigenvalue of the matrix $Q=M^{-\frac{1}{2}}KM^{-\frac{1}{2}}$. 
However,  if $\sqrt{\nu}\omega$ is large,  then the SBS method may converge slowly or even fail to converge \cite{SalkuyehMATCOM}.
}

\cred{One may rewrite  the system \eqref{MainSystem15} in the real form
\begin{equation}\label{AxLu1}
\mathcal{K}	\begin{pmatrix}
	r   \\
	s  
	\end{pmatrix}
\equiv
	\begin{pmatrix}
		E  & F^T \\
		F  & -E
	\end{pmatrix}
	\begin{pmatrix}
		r   \\
		s  
	\end{pmatrix}
	=
	\begin{pmatrix}
		t   \\
		0  
	\end{pmatrix},
\end{equation}
where 
\begin{equation}\label{ABEq}
	E =
	\begin{pmatrix}
		M & O  \\
		O & M 
	\end{pmatrix},
	\quad
	F=
	\begin{pmatrix}
		\sqrt{\nu}K &-\sqrt{\nu}\omega M \\
		\sqrt{\nu}\omega M & \sqrt{\nu}K 
	\end{pmatrix},
\end{equation}
and 
\[
r=\begin{pmatrix}
\Re(\bar{y})   \\
\Im(\bar{y}) 
\end{pmatrix},\quad
s=\begin{pmatrix}
\Re(\bar{q})   \\
\Im(\bar{q}) 
\end{pmatrix},\quad
t=\begin{pmatrix}
\Re( {\hat{y}_d})   \\
\Im( {\hat{y}_d})
\end{pmatrix}.
\]
It is noted that if $z$ is a complex vector, then  $\Re(z)$ and $\Im(z)$ denote  the real and imaginary parts of $z$, respectively.
Axelsson and Luk\'{a}\v{s} in \cite{PRESB}  (see also \cite{AxelMMA}) applied the preconditioned square block (PRESB) preconditioner 
\begin{equation}\label{AxLu2}
	\mathcal{C}=
	\begin{pmatrix}
		E+F+F^T &  F^T \\
		F    &  -E
	\end{pmatrix},
\end{equation}
for system of linear equations of the form \eqref{AxLu1}.  They proved that every eigenvalue of the matrix  $\mathcal{C}^{-1}\mathcal{K}$ lies in the interval $[\frac{1}{2},1]$ (see \cite{PRESB,SalkyehBIT,NUMAxel,CMA2013}).  Therefore, the GMRES method is indeed appropriate for solving the preconditioned system.  It turns out that for applying the preconditioner $\mathcal{C}$, we need solving two systems where the corresponding  coefficient matrices are $E+F$ and $E+F^T$. We can also solve these systems using  the GMRES method  incorporated with the PRESB preconditioner. To apply the PRESB preconditioner  for solving these systems, we require to solve two subsystems with the coefficient matrix $S=(1+w\sqrt{\nu})M+\sqrt{\nu}K$.  Since the matrix $S$ is symmetric positive definite, these systems can be solved either exactly by the Cholesky factorization or inexactly by the conjugate gradient (CG) method. It is worth nothing that if  these systems are solved inexactly then we must employ the flexible GMRES (FGMRES) \cite{FGMRES} instead of the GMRES method.}

In this paper, we present a stationary iteration method, namely Alternating Symmetric positive definite and  Scaled symmetric positive semidefinite Splitting (ASSS), for solving the system (\ref{MainSystem15}) and verify its convergence properties and the corresponding induced preconditioner.

We use the following notation in the rest of the paper. For a square matrix $Z$, we denote  $\sigma(Z)$ , $\rho(Z)$ and $\|Z\|_2$ for the spectrum, spectral radius and 2-norm of the matrix, respectively.  The \textsc{Matlab} notation $(x;y)$ is used to denote the vector $(x^T,y^T)^T$. The imaginary unit ($\sqrt{-1}$) is shown by  the letter $i$.  
The matrix $A\in\Bbb{R}^{n\times n}$ is called positive real if $x^TAx>0$, for every  $0\neq x\in\Bbb{R}^n$.

This paper is organized as follows. In Section \ref{Sec2} the ASSS iteration method is introduced and its convergence analysis is given. Inexact version of the ASSS iteration method is presented in Section \ref{Sec3}. Section \ref{Sec4} is devoted to introducing the ASSS preconditioner and its implementation issues. Numerical results are presented in Section \ref{Sec5} to demonstrate the efficiency of the ASSS iteration method and the corresponding preconditioner. Some concluding results are drawn in Section \ref{Sec6}.

\section{The ASSS iteration method}\label{Sec2}

Using the idea of  \cite{RSS}, the system \eqref{MainSystem15} can be written in the  4-by-4 block real system
\begin{equation}\label{EqASSS1}
\mathcal{A} {\bf x} \equiv\left(\begin{array}{ccccc}
M            &        0        & \sqrt{\nu}K   &  \omega\sqrt{\nu} M \\
0            &        M        & -\omega\sqrt{\nu}M &  \sqrt{\nu}K   \\
\sqrt{\nu}K  & -\omega\sqrt{\nu}M   &      -M       &      0         \\
\omega\sqrt{\nu}M &    \sqrt{\nu}K  &       0       &    -M   \\ 
\end{array}\right)\left(\begin{array}{l}
\Re(\bar{y}) \\
\Im(\bar{y}) \\
\Re(\bar{q}) \\
\Im(\bar{q}) \\
\end{array}\right)=\left(\begin{array}{c}
\Re(\hat{y}_d) \\
\Im(\hat{y}_d) \\
0\\
0
\end{array}\right) \equiv \hat{\mathbf{b}}.
\end{equation}
We define the matrices ${\mathcal G}_1$ and ${\mathcal G}_2$ as following
\[
{\mathcal G}_1=\left(\begin{array}{ccccc}
I            &        0        &           0   &  \omega\sqrt{\nu} I \\
0            &        I        & -\omega\sqrt{\nu}I &  0   \\
0            & -\omega\sqrt{\nu}I   &      -I       &      0         \\
\omega\sqrt{\nu}I &        0        &       0       &    -I   \\ 
\end{array}\right),\quad 
{\mathcal G}_2=\sqrt{1+\nu \omega^2}\left(\begin{array}{ccccc}
0   &  0  &  I   &  0  \\
0   &  0  &  0   &  I  \\
I   &  0  &  0   &  0  \\
0   &  I  &  0   &  0  \\ 
\end{array}\right),\quad 
\]
where $I\in\Bbb{R}^{m\times m}$ is the identity matrix. Then, the system \eqref{EqASSS1} can be equivalently rewritten as
\begin{equation}\label{EqASSS2}
({\mathcal G}_1{\mathcal M}+\frac{\sqrt{\nu}}{\sqrt{1+\nu \omega^2}} {\mathcal G}_2 \hat{{\mathcal K}} ){\bf x}
=\hat{{\bf b}}.
\end{equation}
where 
\[
{\mathcal M}=\left(\begin{array}{cccc}
M  &  0   &  0   &  0   \\
0  &  M   &  0   &  0   \\
0  &  0   &  M   &  0   \\
0  &  0   &  0   &  M  \\ 
\end{array}\right)\quad  
{\rm and} \quad
\hat{{\mathcal K}}=\left(\begin{array}{cccc}
K  &  0   &  0   &  0   \\
0  &  K   &  0   &  0   \\
0  &  0   &  K   &  0   \\
0  &  0   &  0   &  K  \\ 
\end{array}\right).
\]
\begin{lemma}
The matrix  ${\mathcal G}_1$ is nonsingular and 
\[
{\mathcal G}_1^{-1}=\frac{1}{1+\nu \omega^2} 
\left(\begin{array}{ccccc}
I                   &  0                 &  0                    &  \omega \sqrt{\nu}I  \\
0                   &  I                 &  -\omega \sqrt{\nu}I   &  0  \\
0                   & -\omega \sqrt{\nu}I &  -I                   &  0    \\
\omega \sqrt{\nu}I   &  0                 &  0                    &  -I  \\ 
\end{array}\right).
\] 
Moreover, the matrix ${\mathcal G}:={\mathcal G}_1^{-1}{\mathcal G}_2$ is of the form 
\[
{\mathcal G}=\frac{1}{\sqrt{\nu(1+\nu \omega^2)}} 
\left(\begin{array}{ccccc}
0                   &  \omega \nu I    &  \sqrt{\nu} I      &  0 \\
-\omega \nu I       &  0                 &  0   &  \sqrt{\nu} I  \\
-\sqrt{\nu} I       & 0  &  0             &  -\omega \nu I    \\
0   &  -\sqrt{\nu}I                 &  \omega \nu I                     &  0  \\ 
\end{array}\right).
\] 
Furthermore, we have ${\mathcal G}^2=-I$ with $I$ being the identity matrix of order $4m$, ${\mathcal G}^{-1}=-{\mathcal G}$ and ${\mathcal G}^T=-{\mathcal G}$.  
\begin{proof}
	The proof is straightforward and is omitted here. 
\end{proof}
\end{lemma}

Premultiplying both sides of the system \eqref{EqASSS2} by ${\mathcal G}_1^{-1}$, gives the system 
\begin{equation}\label{EqASSS3}
\mathcal{B}\bf{x}\equiv ({\mathcal M}+ {\mathcal G} {\mathcal K}){\bf x}
={\bf b},
\end{equation}
where  $\mathcal{K}=\eta \hat{\mathcal{K}}$ with $\eta=\sqrt{\nu}/\sqrt{1+\nu \omega^2}$ and ${\bf b}={\mathcal G}_1^{-1}\hat{\bf{b}}$.

\cblu{Given $\alpha>0$, we rewrite Eq. \eqref{EqASSS3} as}
\begin{equation}\label{EqASSSA}
(\alpha I+{\mathcal M}){\bf x}=(\alpha I -{\mathcal G} {\mathcal K}) {\bf x}+{\bf b}.
\end{equation}
\cblu{Eq. \eqref{EqASSS3} can also be reformulated in the following form}
\[
\mathcal{G} (\alpha I  +  {\mathcal K}){\bf x}=(\alpha {\mathcal G} -\mathcal{M}) {\bf x}+{\bf b}.
\]
Premultiplying both sides of this equation by $\mathcal{G}^{-1}$ and having in mind that $\mathcal{G}^{-1}=-\mathcal{G}$, gives
\begin{equation}\label{EqASSSB}
(\alpha I  +  {\mathcal K}){\bf x}=(\alpha I +\mathcal{G}\mathcal{M}) {\bf x}-\G{\bf b}.
\end{equation}
Now, using Eqs. \eqref{EqASSSA} and \eqref{EqASSSB} we establish the ASSS iteration method for solving the system \eqref{EqASSS3} as
\begin{equation}\label{ASSS}
\left\{
\begin{array}{rl}
(\alpha I+{\mathcal M}){\bf x}^{(k+\frac{1}{2})}\hspace{-0.25cm} &=(\alpha I -{\mathcal G} {\mathcal K}) {\bf x}^{(k)}+{\bf b},\\
(\alpha I  +  {\mathcal K}){\bf x}^{(k+1)}\hspace{-0.25cm}&=(\alpha I +\mathcal{G}\mathcal{M}) {\bf x}^{(k+\frac{1}{2})}-\G{\bf b},
\end{array}
\right.
\end{equation}
where $\bf{x}^{(0)}$ is an initial guess. \cblu{It is worth noting that  both of the matrices ${\mathcal M}$ and ${\mathcal K}$ are  
	symmetric positive definite and the ASSS method is obtained using the splitting $\mathcal{B}={\mathcal M}+ {\mathcal G} {\mathcal K}$, which is indeed a splitting of a symmetric positive definite matrix and a scaled (by $\mathcal{G}$) symmetric positive definite matrix.}

  Eliminating ${\bf x}^{(k+\frac{1}{2})}$ from Eq. \eqref{ASSS} yields
the following stationary iterative method
\begin{equation}\label{ASSSSt}
{\bf x}^{(k+1)}=\T_{\alpha} {\bf x}^{(k)}+\bf{f},
\end{equation}
where 
\[
\T_{\alpha}=(\alpha I+\K)^{-1} (\alpha I+\G\M) (\alpha I+\M)^{-1} (\alpha I-\G\K),
\] 
is the iteration matrix of the ASSS method and $${\bf f}=\alpha(\alpha I+\K)^{-1}(I-\G)(\alpha I+\M)^{-1}{\bf b}.$$  The following theorem states the convergence of the ASSS method. 
\begin{theorem}\label{Thm1}
	Assume that the matrices $K$ and $M$ are  symmetric positive definite, and let $\alpha>0$. Then the spectral radius $\rho(\T_{\alpha})$ of the ASSS iteration matrix 
	\[
	\T_{\alpha}=(\alpha I+\K)^{-1} (\alpha I+\G\M) (\alpha I+\M)^{-1} (\alpha I-\G\K),
	\]  
	satisfies 
	\[
	\rho(\T_\alpha)\leq \gamma(\alpha)=\max_{\mu \in \sigma(M)} \frac{\sqrt{\alpha^2+\mu^2}}{\alpha+\mu}\max_{\lambda \in \sigma(K)} \frac{\sqrt{\alpha^2+(\eta\lambda)^2}}{\alpha+\eta\lambda},
	\] 
	where $\eta=\sqrt{\nu}/\sqrt{1+\nu \omega^2}$.	Hence, it holds that
	\[
	\rho(\T_\alpha)\leq \gamma(\alpha)< 1,\qquad \forall \alpha>0,
	\]
	which shows that the ASSS iteration method converges unconditionally.
	\begin{proof}
		By setting $\hat{\T}_{\alpha} = (\alpha I+\K) \T_{\alpha} (\alpha I+\K)^{-1}$, we see that
		\begin{eqnarray*}
		\hat{\T}_{\alpha} &=& (\alpha I+\G\M) (\alpha I+\M)^{-1} (\alpha I-\G\K) (\alpha I+\K)^{-1} \\
		                  &=& \R_{\alpha}\S_{\alpha},
		\end{eqnarray*}
		where	$\R_{\alpha}=(\alpha I+\G\M)(\alpha I+\M)^{-1}$ and $\S_{\alpha}=(\alpha I-\G \K)(\alpha I+\K)^{-1}$. Note that
		the matrices $\hat{\T}_{\alpha}$ and ${\T}_{\alpha}$ are similar, and as a result we deduce that	 
		\[
		\rho(\T_{\alpha})=\rho(\hat{\T}_{\alpha})=\rho(\R_{\alpha} \S_{\alpha})\leq \|\R_{\alpha} \S_{\alpha}\|_2\leq \|\R_{\alpha}\|_2 \|\S_{\alpha}\|_2.
		\]
		On the other hand, it follows from $\G^2=-I$, $\G^T=-\G$ and $\M\G=\G\M$ that 
		\begin{eqnarray}
		\nonumber\|\R_{\alpha}\|_2^2 &=& \rho \left((\alpha I+\M)^{-1} (\alpha I-\M\G) (\alpha I+\G\M) (\alpha I+\M)^{-1}\right)\\
		\nonumber                    &=& \rho ((\alpha I+\M)^{-1} (\alpha^2 I+\alpha\M\G -\alpha \G \M-\M\G^2 \M) (\alpha I+\M)^{-1})\\
		\nonumber                    &=& \rho \left((\alpha I+\M)^{-2} (\alpha^2 I+\M^2) (\alpha I+\M)^{-1}\right)\\
		\nonumber &=& \max_{\mu \in \sigma(\M)} \frac{\alpha^2+\mu^2}{(\alpha+\mu)^2}\\
		&=&\max_{\mu \in \sigma(M)} \frac{\alpha^2+\mu^2}{(\alpha+\mu)^2}. \label{Up1}
		\end{eqnarray}
		In the same way, we deduce that
		\begin{equation}\label{Up2}
		\|\S_{\alpha}\|_2^2 = \max_{\lambda \in \sigma(K)} \frac{\alpha^2+(\eta\lambda)^2}{(\alpha+\eta\lambda)^2}.
		\end{equation}
		Since the matrices $M$ and $K$ are symmetric positive definite we have $\mu,\lambda>0$, for all $\mu\in\sigma(M)$ and $\mu\in\sigma(K)$, and as result we deduce that $\sqrt{\alpha^2+\mu^2}< \alpha+\mu$ and 
		 $\sqrt{\alpha^2+(\eta\lambda)^2}< \alpha+\eta\lambda$. Therefore, we see that $\|\R_{\alpha}\|_2<1$ and $\|\S_{\alpha}\|_2<1$. Hence,  $\rho(\T_\alpha)\leq \gamma({\alpha})=\|\R_{\alpha}\|_2\|\S_{\alpha}\|_2 <1$, which completes the proof.
	\end{proof}
\end{theorem}

\begin{remark} 
Under the conditions of Theorem \ref{Thm1}, we have
	\[
	\max_{\lambda \in \sigma(K)} \frac{\sqrt{\alpha^2+(\eta\lambda)^2}}{\alpha+\eta\lambda}<1,\quad \forall \alpha>0.
	\] 
Hence,	
	\begin{equation}\label{Zeta}
	\rho(\T_\alpha)\leq \zeta(\alpha)=\max_{\mu \in \sigma(M)} \frac{\sqrt{\alpha^2+\mu^2}}{\alpha+\mu}.
	\end{equation}
Similar to Corollary 2.1 in \cite{MHSS} the minimum value of the $\zeta({\alpha})$ is obtained at 
\begin{equation}\label{alphast}
\alpha^*=\sqrt{\mu_{\min}\mu_{\max}},
\end{equation}
where $\mu_{\min}$ and $\mu_{\max}$ are the smallest and largest eigenvalues of the matrix $M$, respectively.
It is worth noting that the parameter $\alpha^*$ is independent of the parameters $\nu$ and $\omega$. 
\cred{In practice, we can use a few iterations of the power method for computing $\mu_{\min}$ and  the inverse power method to calculate  $\mu_{\max}$.}
\end{remark}

\begin{remark}
	In Theorem  \ref{Thm1} if the matrix $K$ is symmetric positive semidefinite, then
		\[
		\max_{\lambda \in \sigma(K)} \frac{\sqrt{\alpha^2+(\eta\lambda)^2}}{\alpha+\eta\lambda}\leq 1,\quad \forall \alpha>0,
		\] 
and  we have 	$\rho(\T_\alpha)\leq \zeta(\alpha)<1$, where  $\zeta(\alpha)$ was defined in Eq. \eqref{Zeta}.
Hence, the ASSS iteration method is still convergent unconditionally.
\end{remark}

\cred{There are several methods to estimate the parameters of the iterative methods and their induced preconditioners having similar structure to the ASSS iteration method (see \cite{Huang2014225,Chen2015271,Yang201825,Ren201636,Cao2020374}). In \cite{Ren201636}, Ren and Cao proposed an efficient estimation formula 	for the iteration parameter of a class of alternating positive-semidefinite splitting preconditioner. The same idea was used by Cao to estimate the parameter of the block positive-semidefinite splitting preconditioner for generalized saddle point linear systems \cite{Cao2020374}.	Similar to the strategy used in \cite{Ren201636,Cao2020374}, we can estimate the parameter of the ASSS method. However, 	it is not as effective as the parameter $\alpha^*=\sqrt{\mu_{\min}\mu_{\max}}$. We will shortly see in Section \ref{Sec5} that the value of $\alpha^*$ can be computed inexpensively.}

\section{Inexact version of ASSS}\label{Sec3}

Let ${\bf x}^{(k)}$ be the computed solution solution at iteration $k$. Setting 
\[
{\bf x}^{(k+\frac{1}{2})}={\bf x}^{(k)}+{\delta}^{(k)},
\]
and substituting it in the first step of Eq. \eqref{ASSS} gives 
\[
(\alpha I+{\mathcal M}){\delta}^{(k)}={\bf b}-{\cal B}{\bf x}^{(k)}=:{\bf r}^{(k)}.
\]
In the same way if we set ${\bf x}^{(k+1)}={\bf x}^{(k+\frac{1}{2})}+{\delta}^{(k+\frac{1}{2})}$ and substitute it in the second step of Eq. \eqref{ASSS}, then we get
\begin{eqnarray*}
\cblu{(\alpha I  +  {\mathcal K}){\delta}^{(k+\frac{1}{2})}} & =&  -{\mathcal K} {\bf x}^{(k+\frac{1}{2})}+\mathcal{G}\mathcal{M} {\bf x}^{(k+\frac{1}{2})}-\G{\bf b} \\
& = &  {\mathcal G}^2 {\mathcal K} {\bf x}^{(k+\frac{1}{2})}+\mathcal{G}\mathcal{M} {\bf x}^{(k+\frac{1}{2})}-\G{\bf b} \\
& = &  -\G \left( {\bf b}- (\M+\G\K) {\bf x}^{(k+\frac{1}{2})} \right)\\
& = &  -\G \left( {\bf b}- \B {\bf x}^{(k+\frac{1}{2})} \right)\\
& = &  -\G {\bf r}^{(k+\frac{1}{2})}.
\end{eqnarray*} 

Now using the above results we establish the inexact version of ASSS (IASSS) as following.

\medskip

\indent {\bf Algorithm 1. The IASSS algorithm for solving $\B{\bf x}={\bf b}$.}\vspace{-0.2cm}
\begin{enumerate}
	\item Choose an initial guess ${\bf x}^{(0)}$. \\[-0.68cm]
	\item For $k=0,1,2,\ldots$ until convergence, Do\\[-0.68cm]
	\item \qquad Compute ${\bf r}^{(k)}={\bf b}-\B {\bf x}^{(k)}$.\\[-0.68cm]
	\item \qquad Solve $(\alpha I+{\mathcal M}){\delta}^{(k)}={\bf r}^{(k)}$ approximately using an iteration method.\\[-0.68cm]
	\item \qquad  ${\bf x}^{(k+\frac{1}{2})}:={\bf x}^{(k)}+\delta^{(k)}$.\\[-0.68cm]
	\item \qquad Compute ${\bf r}^{(k+\frac{1}{2})}={\bf b}-\B {\bf x}^{(k+\frac{1}{2})}$.\\[-0.68cm]
	\item \qquad Solve $(\alpha I+\K){\delta}^{(k+\frac{1}{2})}=-\G{\bf r}^{(k+\frac{1}{2})}$ approximately using an iteration method.\\[-0.68cm]
	\item \qquad  ${\bf x}^{(k+1)}:={\bf x}^{(k+\frac{1}{2})}+\delta^{(k+\frac{1}{2})}$.\\[-0.68cm]
	\item EndDo
\end{enumerate}

In the steps 4 and 7 of Algorithm 1 two \cblu{linear systems of equations} with the coefficient matrices  $\alpha I+{\mathcal M}$ and $\alpha I+{\mathcal K}$ should be solved. Since these matrices are SPD, they can be solved inexactly using the CG method. However, since 
\[
\alpha I+{\mathcal M}=
\left(\begin{array}{cccc}
	\alpha I+ M  &  0   &  0   &  0   \\
	0  &  \alpha I+ M   &  0   &  0   \\
	0  &  0   &  \alpha I+ M   &  0   \\
	0  &  0   &  0   &  \alpha I+ M  \\ 
\end{array}\right),
\] 
we need to solve four systems with the same coefficient matrices $\alpha I+ M$ and different right-hand sides. In fact if we set \cblu{${\hat\delta}^{(k)}=[{\delta}^{(k)}_1;{\delta}^{(k)}_2;{\delta}^{(k)}_3;{\delta}^{(k)}_4]$ $\in\Bbb{R}^{4m}$ 
and ${\hat{\bf r}}^{(k)}=[{\bf r}^{(k)}_1;{\bf r}^{(k)}_2;{\bf r}^{(k)}_3;{\bf r}^{(k)}_4]\in\Bbb{R}^{4m}$} with 
${\delta}^{(k)}_i,{\bf r}^{(k)}_i\in\Bbb{R}^m$, for $i=1,2,3,4$, then we only need to solve the linear 
system of equations with multiple right-hand sides 
\[
(\alpha I+ M){\hat\delta}^{(k)}={\hat{\bf r}}^{(k)}.
\]
This system can be solved using the global CG algorithm with a suitable preconditioner, which can be written as following. In this algorithm the inner product used is $\langle X,Y \rangle=trace(X^TY)$.

\medskip

\indent {\bf Algorithm 2.  Global CG algorithm for $P^{-1}AX=P^{-1}B$.}\vspace{-0.2cm}
\begin{enumerate}
	\item Compute $R_0=B-AX_0$, $Z_0=P^{-1}R_0$ and $P_0=Z_0$. \vspace{-0.2cm}
	
	\item For $j:=0,1,\ldots,$ until convergence, Do\vspace{-0.2cm}
	
	\item \qquad $\alpha_j:=\langle R_j,Z_j\rangle/\langle AP_j,P_j\rangle $.
	\vspace{-0.2cm}
	
	\item \qquad $X_{j+1}:=X_j+\alpha_jP_j$.\vspace{-0.2cm}
	
	\item \qquad $R_{j+1}:=R_j-\alpha_jAP_j$.  \vspace{-0.2cm}
	
	\item \qquad $Z_{j+1}:=P^{-1}R_{j+1}$.\vspace{-0.2cm}
		
	\item \qquad $\beta_j:=\langle R_{j+1},Z_{j+1} \rangle/ \langle R_j,Z_j\rangle$.
	\vspace{-0.2cm}
	
	\item \qquad $P_{j+1}:=Z_{j+1}+\beta_jP_j$. \vspace{-0.2cm}
	
	\item Enddo
\end{enumerate}

In the same way the  system of Step 7 of Algorithm 1 can be solved. So in each iteration of the IASSS iteration method two linear systems of equations with multiple right-hand sides should be solved which can be solved using Algorithm 2. \cblu{The incomplete Cholesky factorization of the coefficient matrix can be used as the preconditioner in Algorithm 2.}
  
As Benzi and Golub mentioned in \cite{Benzi-SD}, adding $\alpha>0$ to the main diagonal of 
$M$ (or $\eta K$)  improves the condition number of the matrix. So this improves, in turn, the  convergence of the CG method much. 

\section{The ASSS preconditioner}\label{Sec4}

It follows from $\G^T=-\G$ and $\G\K=\K\G$ that
\[
\B^T+\B=\M+ \G \K+\M+\K\G^T=\M+ \G \K+\M-\K\G=2\M,
\]
which shows that the matrix $\B$ is positive definite (positive real). Hence, from Theorem 6.30 in \cite{SaadBook} the restarted  GMRES($\ell$) \cblu{for solving the system \eqref{EqASSS3}} converges for any $\ell\geq 1$.    
On the other hand, if we define 
\begin{eqnarray*}
\P_{\alpha}&=&\frac{1}{\alpha} (I+\G)^{-1} (\alpha I+\M) \G (\alpha I+\K),\\
\Q_{\alpha}&=&\frac{1}{\alpha} ( I+\G)^{-1} (\alpha \G-\M) (\alpha I-\G\K),
\end{eqnarray*}
then 
\[
\B=\P_{\alpha}-\Q_{\alpha}\quad \text{and} \quad \T_{\alpha}=\P^{-1}_{\alpha}\Q_{\alpha}.
\]
 Hence, we get
\[
\P^{-1}_{\alpha}\B=I-\T_{\alpha},
\]
which shows that the eigenvalues of $\P^{-1}_{\alpha}\B$ are clustered in a circle with radius 1 centered at $(1,0)$.
Hence, the GMRES method would be quite appropriate for solving the system incorporated with the preconditioner $\P_{\alpha}$, i.e., 
\[
\P_{\alpha}^{-1}\B {\bf x}
=\P_{\alpha}^{-1}{\bf b}.
\]

In applying the preconditioner in each iteration of a Krylov subspace method like GMRES a linear system of equations of the form $\P_{\alpha}s=r$ should be solved. Since 
\begin{eqnarray*}
\P_{\alpha}^{-1}&=&\alpha (\alpha I+\K)^{-1} \G^{-1}  (\alpha I+\M)^{-1} (I+\G)\\
&=&-\alpha (\alpha I+\K)^{-1} \G  (\alpha I+\M)^{-1} (I+\G),
\end{eqnarray*}  
we can state the following algorithm for solving the system $\P_{\alpha}s=r$.

\bigskip

\indent {\bf Algorithm 3. Solution of $\P_{\alpha}s=r$.}\vspace{-0.2cm}
\begin{enumerate}
	\item Compute $v=-\alpha(I+\G)r$. \vspace{-0.2cm}
	
	\item Solve \cblu{$(\alpha I+\M)w=v$ for $w$}.\vspace{-0.2cm}
	
	\item Compute $z=\G w$.\vspace{-0.2cm}
	
	\item Solve \cblu{$(\alpha I+\K)s=z$ for $s$}.
\end{enumerate}

Both of the systems in Steps 2 and 4 of the  above algorithm can be solved exactly using the Cholesky factorization of the matrices $\alpha I+M$ and  $\alpha I+\eta K$, or inexactly using the CG method. When these systems are solved inexactly  we can apply the global CG algorithm described in the previous section.  It is noted that, in this case the flexible version of GMRES (FGMRES) \cite{FGMRES} should be used for solving the main system instead of the GMRES algorithm.

\section{Numerical results}\label{Sec5}

We consider the distributed control problem introduced in Section \ref{Sec1} in two-dimensional case with the computational domain $\Omega=(0,1)\times (0,1)\in \Bbb{R}^2 $. The target state is set to be
\begin{equation}\label{Eq00108}
y_d(x,y)=
\left\{
\begin{array}{cl}
(2x-1)^2(2y-1)^2,&\text{if}~~(x,y)\in (0,\frac{1}{2})\times (0,\frac{1}{2}), \\
0,&\text{otherwise}.
\end{array} \right.
\end{equation}
\cred{In our numerical test, we discretize the problem using the bilinear quadrilateral {\bf Q1} finite elements with a uniform mesh \cite{Elman}.
Let $D=\diag(M)$, where $M$ is the mass matrix. In this case, the smallest and largest eigenvalues of $N=D^{-1}M$ are $\frac{1}{4}$ and $\frac{9}{4}$, respectively \cite{Wathen19877}. On the other hand, $D$ is a scalar multiplication of an identity matrix, i.e., $D=\vartheta I$ for some $\vartheta>0$. Therefore, we have  $M=\vartheta N$. So, we deduce that the smallest and largest eigenvalues of $M$ are $\mu_{\min}=\frac{1}{4}\vartheta$ 
and $\mu_{\max}=\frac{9}{4}\vartheta$, respectively. Therefore, we have
\[
\alpha^*=\sqrt{\mu_{\min} \mu_{\max}}=\frac{3}{4}\vartheta.
\]
It is worth noting that the parameter $\alpha^*$ depends only on the mesh size and independent of the parameters $\nu$ and $\omega$. In Table \ref{Tbl0}, we disclose the values of $\vartheta$, $\mu_{\min}$, $\mu_{\max}$ and $\alpha$ for $h=2^{-k}$, $k=4,5,6,7$.

\begin{table}[!h]
	\centering\caption{The values of $\vartheta$, $\mu_{\min}$, $\mu_{\max}$ and $\alpha$ for $h=2^{-k}$, $k=4,5,6,7$. \label{Tbl0}}
	\begin{center}
		\label{exact}
		\scalebox{1}
		{
		\begin{tabular}{|c|c|c|c|c|} \hline
		$h$         &    $\vartheta   $  & $\mu_{\min}$ &  $\mu_{\max}$ &  $\alpha^*$    \\ \hline
		$2^{-4}$	&    $1.7361\times 10^{-3}$ &  $4.3403\times 10^{-4}$ & $3.9063\times 10^{-3}$ &  $1.3021\times 10^{-3}$\\	
		$2^{-5}$	& ~  $4.3403\times 10^{-4}$  ~  & ~ $1.0851\times 10^{-4}$ ~  &~ $9.7656\times 10^{-4}$ ~  & ~ $3.2552\times 10^{-4}$ ~\\	
		$2^{-6}$	&   $1.0851\times 10^{-4}$    &  $2.7127\times 10^{-5}$   & $2.4414\times 10^{-4}$   &  $8.1380\times 10^{-5}$ \\
		$2^{-7}$	&   $2.7127\times 10^{-5}$    &  $6.7817\times 10^{-6}$   & $6.1035\times 10^{-5}$   &  $2.0345\times 10^{-5}$ \\ \hline
			\end{tabular}}
		\end{center}
	\end{table}}

To generate the system \eqref{MainSystem} we have used the codes of the paper \cite{Rees} which  is
available at \url{www.numerical.rl.ac.uk/people/rees/}. \cred{All runs are implemented in \textsc{Matlab} R2017, equipped with a Laptop with 1.80 GHz central processing unit (Intel(R) Core(TM) i7-4500), 6 GB RAM and Windows 7 operating system.} \cblu{For each method, we report the number of iterations for the convergence and the elapsed CPU time (in seconds). In the tables, a dagger ($\dag$) and a double dagger  ($\ddag$)  mean that the method has not converged in 500 iterations and 60 seconds, respectively. }

We divide the numerical results in two parts. In the first part, we compare the numerical results of the inexact version of the ASSS iteration method described in Section \ref{Sec3} (denoted by IASSS) with those of the inexact version of the BAS iteration method (denoted by IBAS). It is noted that in each iteration of the BAS iteration method two subsystems with the coefficient matrix $\alpha M+\sqrt{\nu} K$ and two systems with the coefficient matrix $M$ should be solved. So in the two half-steps of the IBAS iteration method the subsystems are solved using the global CG algorithm. The outer iteration is terminated as soon as the residual norm of the system \eqref{MainSystem15} is   
reduced by a factor of $10^{6}$. The global CG iteration for solving the subsystems are stopped as soon as the residual Frobenious norm of the residual matrix is reduced by a factor of $10^4$. We always use a zero vector as an initial guess and the maximum number of iterations is set to be 500. 

\cred{In the IASSS iteration method we use the $\alpha^*$ computed by the formula \eqref{alphast} (presented in Table \ref{Tbl0}), and in the IBAS iteration method the parameter $\alpha$ is set to be  $1+\nu \omega^2$ (as suggested in\cite{ZhengCAM}.)} Numerical results for different values of $h$, $\nu$ and $\omega$ have been presented in Tables \ref{Tbl1}, \ref{Tbl2}  and \ref{Tbl3}. 
As we observe there is no significant difference between the number iterations of the IASSS method when the 
parameters $h$, $\nu$ and $\omega$ are changed. Comparing the numerical results of the IASSS method with those of the IBAS method shows that the IBAS iteration method sometimes fails, especially when $\nu$ and $\omega$ are large. However, this is not the case for the IASSS iteration method. Nevertheless, when the values of $\omega$ and $\nu$ are small enough, the CPU time of IBAS is often less than that of IASSS.   
      
\begin{figure}[!h]
	\centering
	\includegraphics[width=8cm,height=6cm]{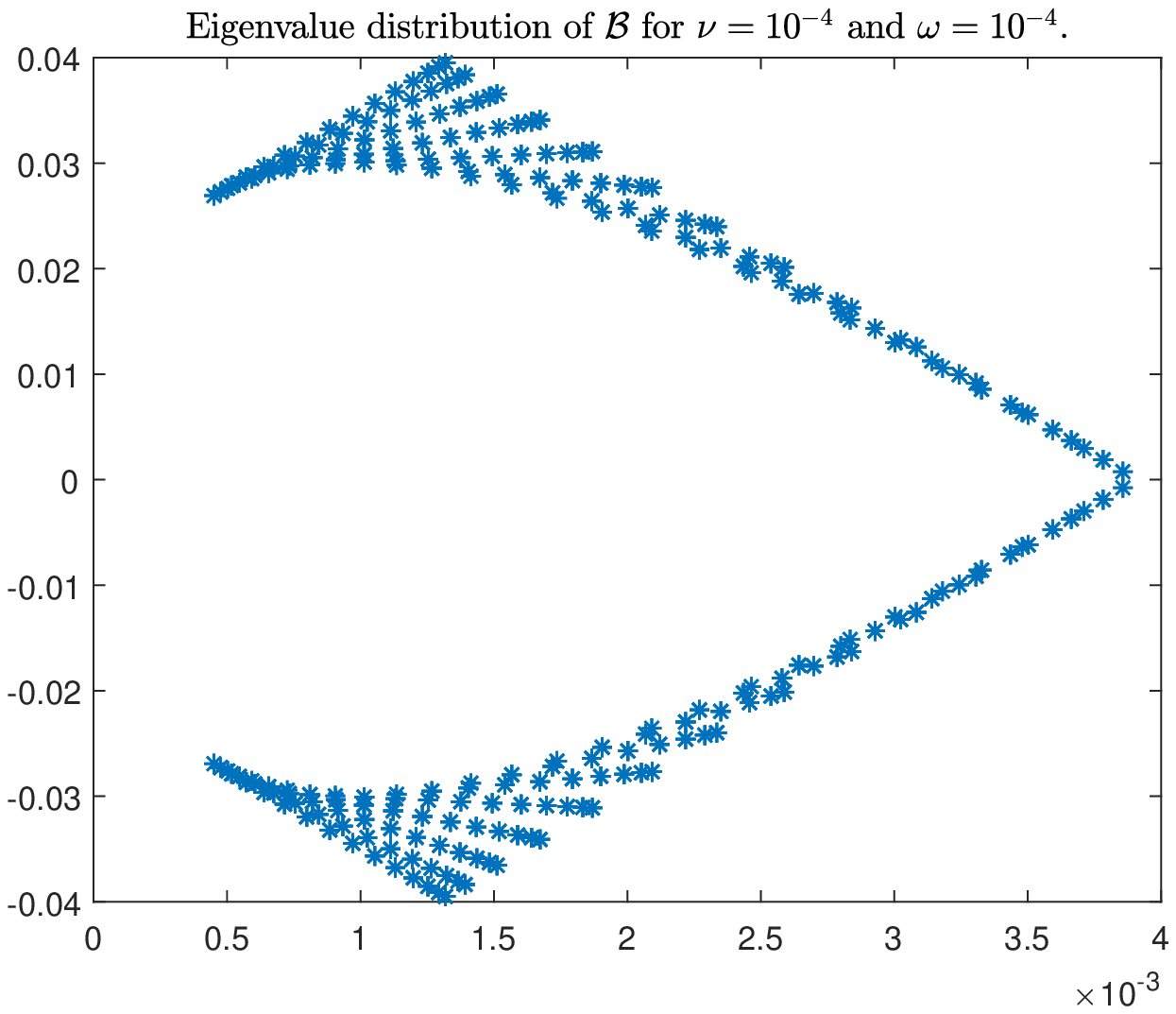}\includegraphics[width=8cm,height=6cm]{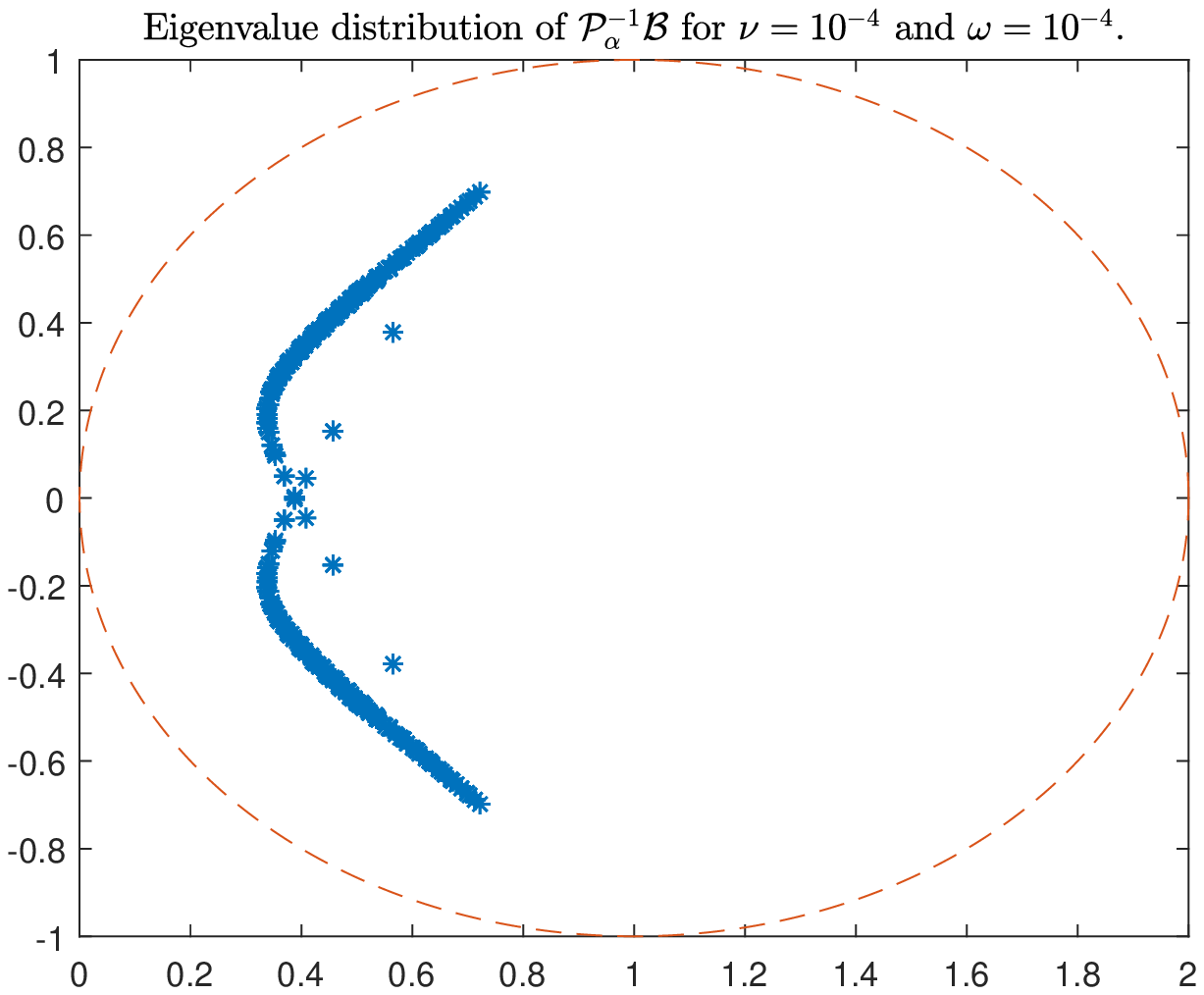}
	\includegraphics[width=8cm,height=6cm]{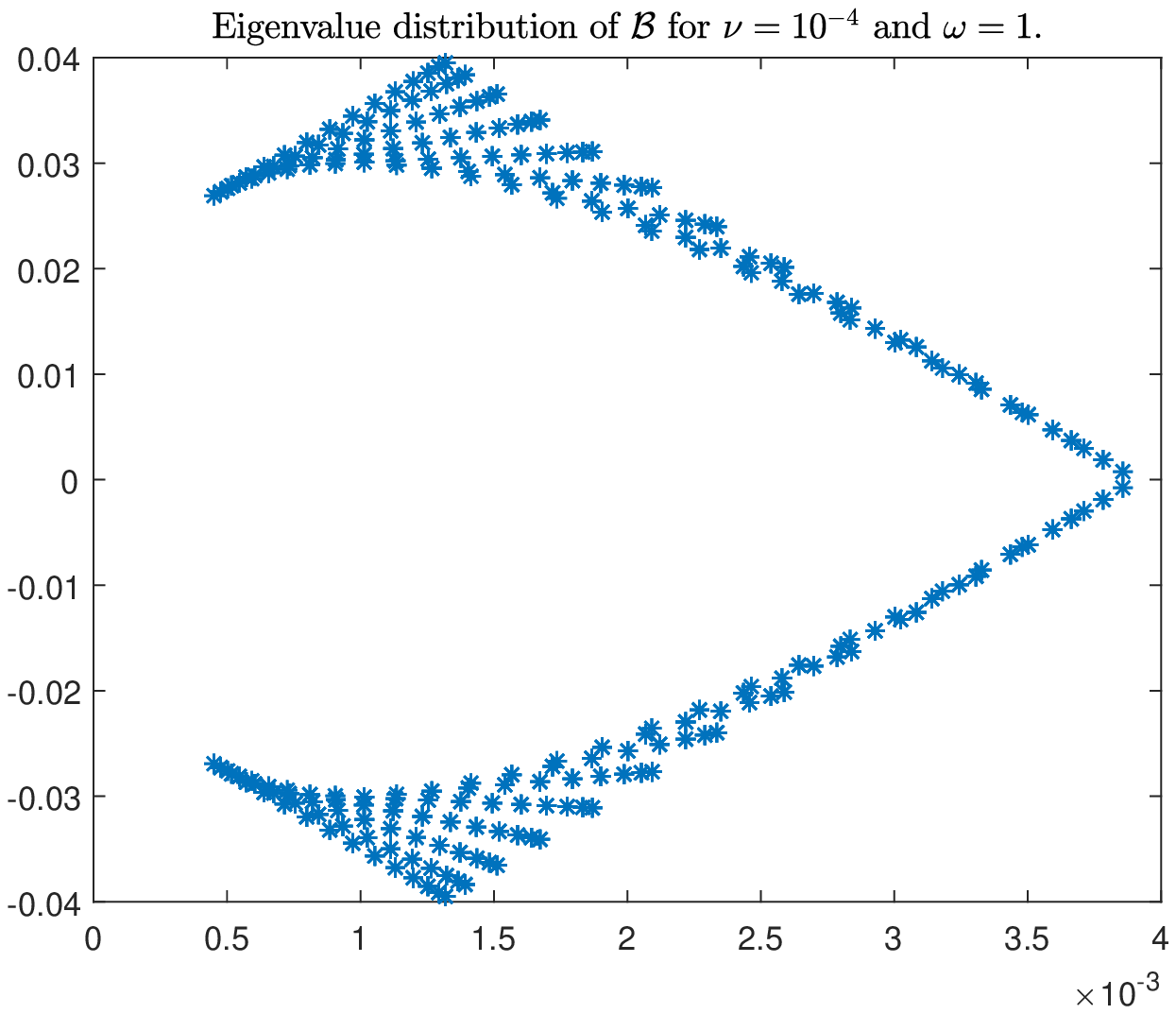}\includegraphics[width=8cm,height=6cm]{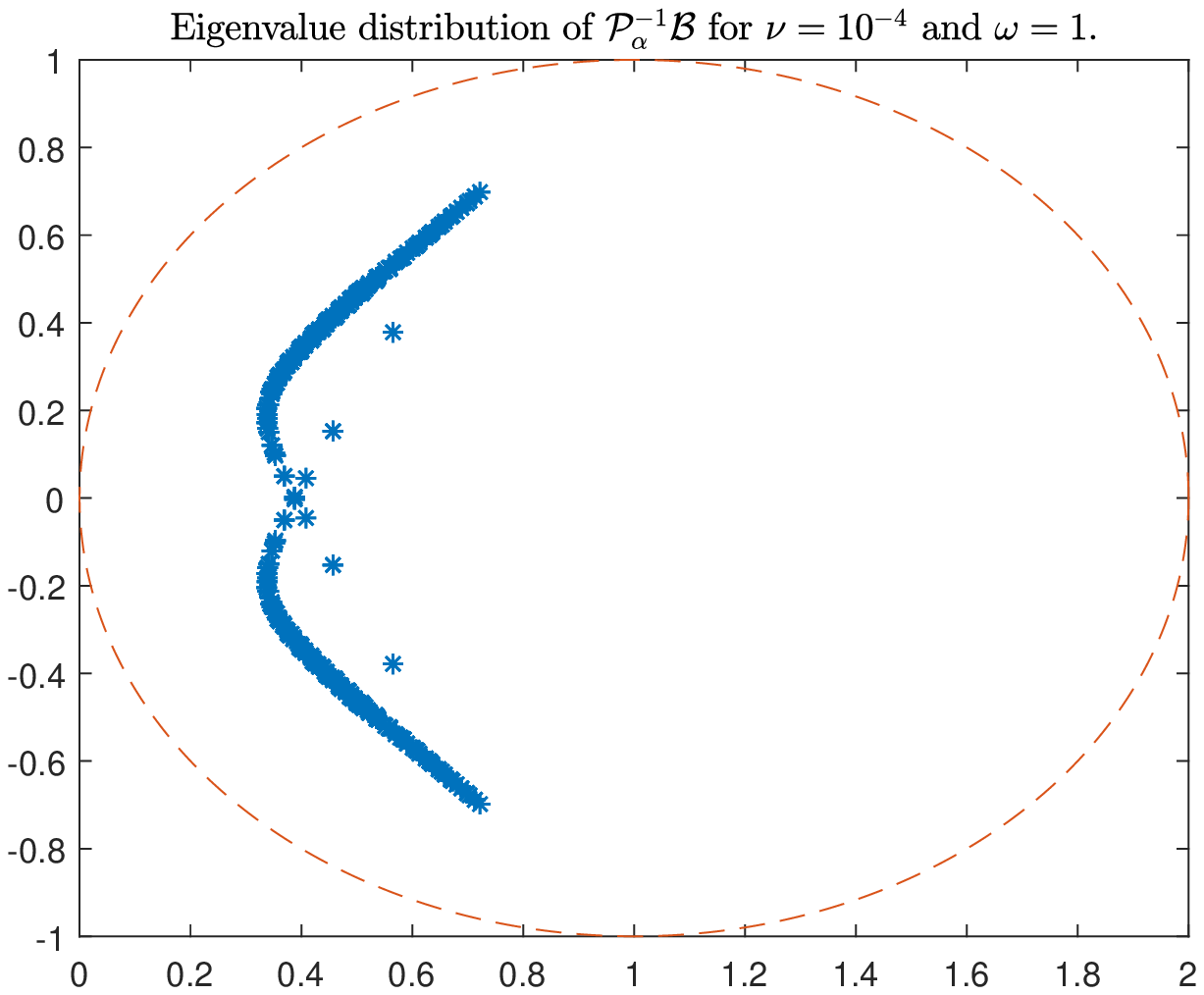}
	\includegraphics[width=8cm,height=6cm]{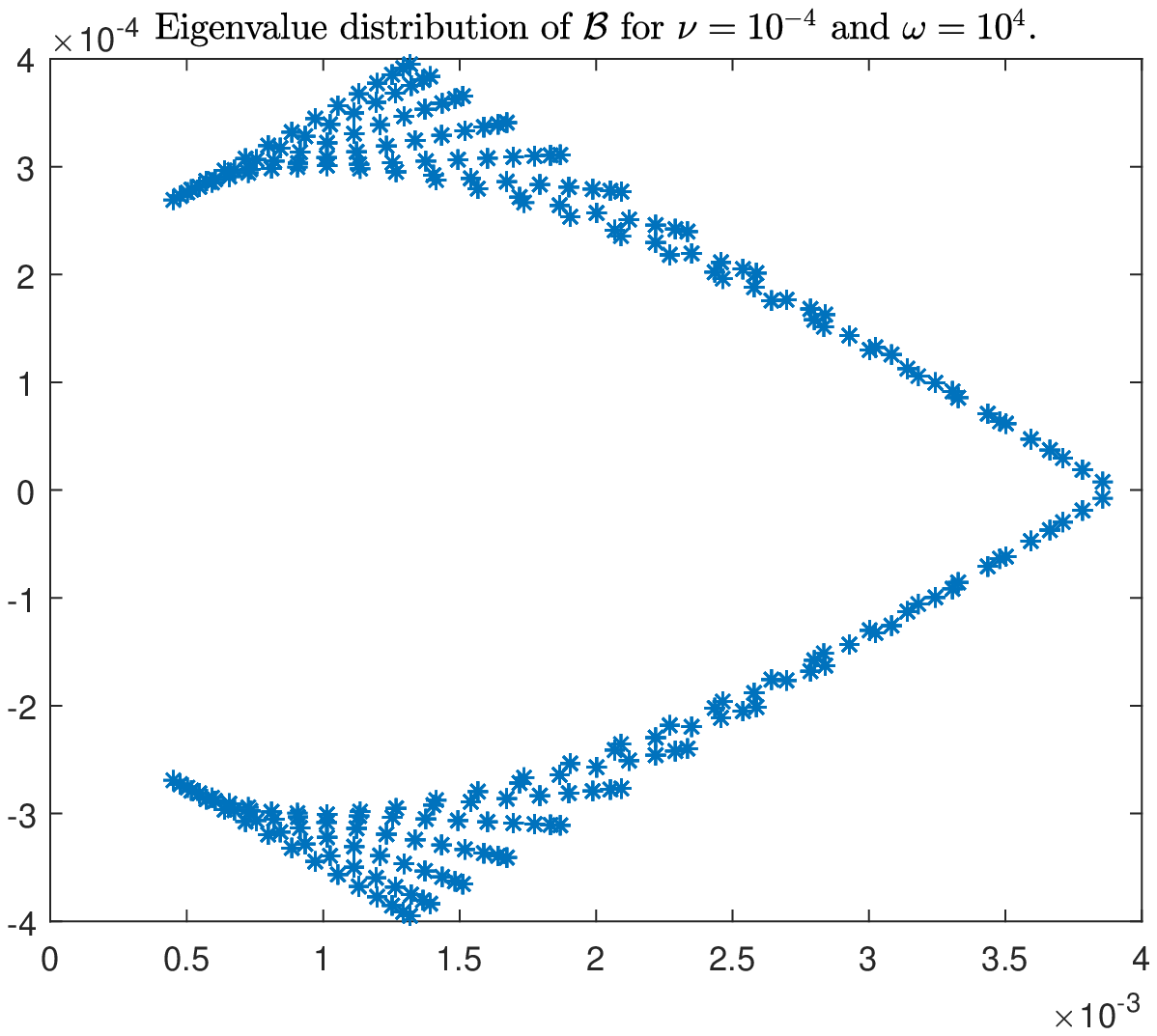}\includegraphics[width=8cm,height=6cm]{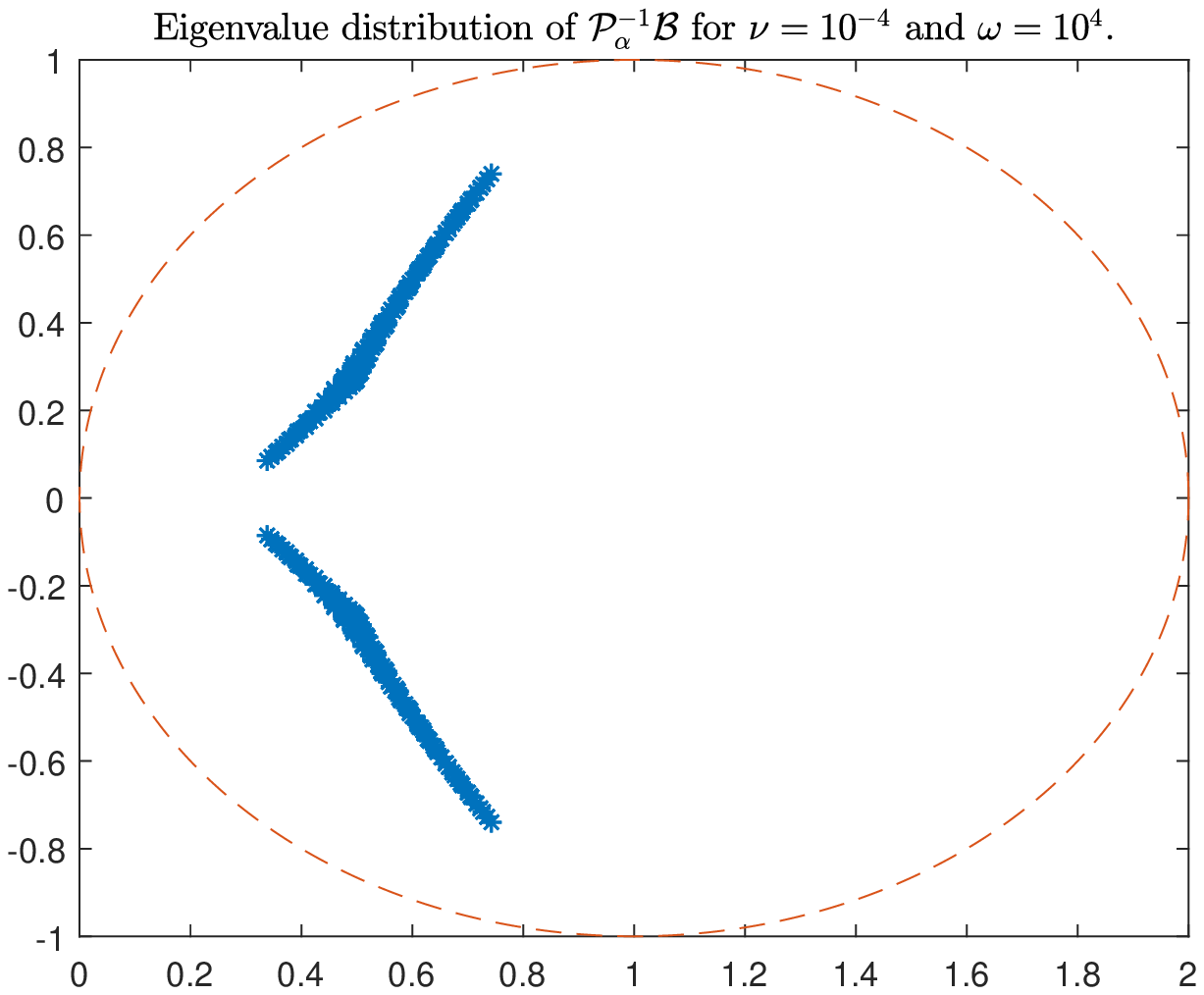}
	\caption{Eigenvalue distribution of the matrices $\mathcal{B}$ and $\mathcal{P}_{\alpha}^{-1}\mathcal{B}$ 
		for $h=2^{-4}$ with different values of $\nu$ and $\omega$. The value of $\alpha$ is $\alpha^*=\sqrt{\mu_{\min}\mu_{\max}},$
		where $\mu_{\min}$ and $\mu_{\max}$ are the smallest and largest eigenvalues of the matrix $M$, respectively. \label{Fig1}}
\end{figure}

\cred{For the second part of our experiments, we compare numerically the performance of the ASSS preconditioner (P-ASSS) with those of the BAS preconditioner (P-BAS), \cred{the PRESB preconditioner (P-PRESB)} and the preconditioner $\mathcal{P}_{BD}$  (P-BD ). To do so, we use the flexible version of the GMRES (FGMRES) method in conjunction with aforementioned preconditioners.
In the implementation of the PRESB preconditioner \eqref{AxLu2}, the systems with the coefficient matrices $E+F$ and $E+F^T$ are solved using the FGMRES method incorporated with the PRESB preconditioner. The innermost subsystems with the coefficient matrix  $S=(1+w\sqrt{\nu})M+\sqrt{\nu}K$ in the PRESB preconditioner and the systems with the coefficient matrix $T=(1+\omega)M+\sqrt{\nu}K$ in $\mathcal{P}_{BD}$
 are solved using the CG method in conjunction with the incomplete Cholesky factorization with dropping tolerance $0.001$ as a preconditioner. For applying the ASSS and the BAS  preconditioners all the subsystems are solved using the global CG method  with incomplete Cholesky factorization with dropping tolerance $0.001$ as a preconditioner. 
We use $\alpha^*$ reported in Table \ref{Tbl0} for the ASSS preconditioner and $\alpha=(1+\nu \omega^2)/(1+\sqrt{\nu}\omega)$ for the BAS preconditioner (as suggested in \cite{ZhengCAM}).

The iteration of the FGMRES method as the outer iteration is stopped as soon as the residual norm is reduced by a factor of $10^6$. All the other iterations are terminated when the residual norm is reduced by a factor of $10^4$. The other  assumptions are as the first part of the numerical experiments.

Numerical results have been presented in Tables \ref{Tbl1}, \ref{Tbl2}  and \ref{Tbl3}.        
As seen, for large values of $\nu$ and $\omega$ the ASSS preconditioner outperforms the other preconditioners. However, for small values of $\nu$ and $\omega$ the ASSS preconditioner is less effective than the others.

Finally, in Figure \ref{Fig1} \cblu{the eigenvalue distribution of the matrices $\mathcal{B}$ and $\mathcal{P}_{\alpha}^{-1}\mathcal{B}$} with $\alpha^*$ (reported in Table \ref{Tbl0}) have been displayed for $h=2^{-4}$ and different values of $\nu$ and $\omega$.
As we observe, the eigenvalues of the preconditioned matrix are well-clustered in the circle with radius 1 and centered at $(1,0)$. }

\begin{table}[!h]
\centering\caption{Number of iterations of the methods along with the elapsed CPU time (in parenthesis)  for $h=2^{-5}$ and different values of $\nu$ and $\omega$. \label{Tbl1}}
\begin{center}
\label{exact}
\scalebox{.8}
{
\begin{tabular}{|c|c|c|c|c|c|c|c|c|c|c|c} \hline
&$\nu\setminus\omega$ &$10^{-4}$ &$10^{-3}$&$10^{-2}$ & $10^{-1}$ & $1$ &$10^{1}$&$10^{2}$&$10^{3}$&$10^{4}$\\ \hline
     
     IASSS& $10^{-2}$& 54(0.12) & 54(0.11) & 54(0.10) & 54(0.10) & 54(0.10) & 53(0.09) & 45(0.07)  & 40(0.06) & 51(0.06) \\
     & $10^{-4}$& 45(0.10) & 45(0.08) & 45(0.07) & 45(0.07) & 45(0.07) & 45(0.07) & 43(0.07)  & 40(0.06) & 51(0.07) \\ 
     & $10^{-6}$& 40(0.09) & 40(0.06) & 40(0.06) & 40(0.06) & 40(0.06) & 40(0.06) & 40(0.06)  & 42(0.06) & 51(0.07) \\
     & $10^{-8}$& 51(0.10) & 51(0.06) & 51(0.07) & 51(0.06) & 51(0.07) & 51(0.06) & 51(0.07)  & 51(0.07) & 52(0.07) \\
\hline
IBAS & $10^{-2}$& 38(0.11) & 38(0.07) & 38(0.10) & 38(0.06) & 38(0.06) & 24(0.04) & 473(0.60) & $\dag$ & $\dag$ \\
     & $10^{-4}$& 33(0.09) & 33(0.05) & 33(0.05) & 33(0.05) & 33(0.06) & 33(0.05) & 39(0.06)  & $\dag$ & $\dag$ \\ 
     & $10^{-6}$& 33(0.08) & 33(0.05) & 33(0.05) & 33(0.06) & 33(0.05) & 33(0.05) & 33(0.05)  & 57(0.07) & $\dag$ \\
     & $10^{-8}$& 38(0.07) & 38(0.06) & 38(0.06) & 38(0.06) & 38(0.06) & 38(0.06) & 38(0.06)  & 38(0.06) & 68(0.08) \\
     \hline\hline\hline
P-ASSS& $10^{-2}$& 29(0.13) & 29(0.09) & 29(0.09) & 29(0.09) & 29(0.09) & 29(0.08) & 28(0.07) & 21(0.05) & 20(0.04) \\
       & $10^{-4}$& 28(0.11) & 28(0.08) & 28(0.07) & 28(0.07) & 28(0.08) & 28(0.07) & 28(0.07) & 21(0.04) & 20(0.04)\\ 
       & $10^{-6}$& 21(0.09) & 21(0.05) & 21(0.05) & 21(0.05) & 21(0.05) & 21(0.05) & 21(0.05) & 21(0.05) & 20(0.04) \\
       & $10^{-8}$& 20(0.09) & 20(0.04) & 20(0.04) & 20(0.04) & 20(0.04) & 20(0.04) & 20(0.04) & 20(0.04) & 19(0.04) \\
\hline
P-BAS & $10^{-2}$& 18(0.12) & 19(0.08) & 20(0.05) & 20(0.05) & 20(0.05) & 17(0.06) & 26(0.07) & 48(0.12) & 36(0.06) \\
& $10^{-4}$& 19(0.12) & 20(0.07) & 21(0.05) & 21(0.05) & 22(0.08) & 21(0.04) & 19(0.04) & 47(0.10) & 40(0.07)\\ 
& $10^{-6}$& 18(0.11) & 19(0.06) & 20(0.04) & 21(0.06) & 21(0.08) & 22(0.04) & 22(0.03) & 28(0.06) & 42(0.08) \\
& $10^{-8}$& 17(0.13) & 18(0.04) & 19(0.03) & 20(0.04) & 21(0.04) & 21(0.04) & 21(0.03) & 21(0.04) & 26(0.04) \\
\hline

P-PRESB& $10^{-2}$& 7(0.20) & 7(0.05) & 7(0.11) & 8(0.10) & 10(0.16) & 23(0.75) & 106(8.80) & $\ddag$  & $\ddag$ \\
       & $10^{-4}$& 8(0.18) & 8(0.04) & 8(0.07) & 8(0.09) &  9(0.09) & 12(0.20) & 60(3.06)  & $\ddag$  & $\ddag$\\ 
       & $10^{-6}$& 8(0.20) & 8(0.05) & 8(0.04) & 8(0.07) &  8(0.06) &  9(0.09) & 12(0.18)  & 89(9.78) & $\ddag$ \\
       & $10^{-8}$& 8(0.19) & 8(0.04) & 8(0.03) & 8(0.04) &  8(0.07) &  8(0.06) & 8(0.05)   & 11(0.17) & 60(3.03) \\
\hline

P-BD   & $10^{-2}$& 14(0.15) & 14(0.07) & 14(0.07) & 16(0.07) & 22(0.08) & 54(0.20) & 176(1.27)  & 132(0.72) & 26(0.06) \\
       & $10^{-4}$& 16(0.14) & 16(0.07) & 18(0.06) & 22(0.09) & 42(0.15) &188(1.42) & 481(25.68) & 351(5.25) & 36(0.07)\\ 
       & $10^{-6}$& 15(0.15) & 15(0.06) & 15(0.04) & 16(0.04) & 30(0.10) &138(0.77) & 345(5.21)  & 298(3.49) & 34(0.07)\\
       & $10^{-8}$& 15(0.13) & 15(0.06) & 15(0.04) & 14(0.03) & 11(0.03) & 22(0.06) & 34(0.08)   & 34(0.07)  & 20(0.04) \\
\hline

\end{tabular}}
\end{center}
\end{table}

\begin{table}[!h]
	\centering\caption{Number of iterations of the methods along with the elapsed CPU time (in parenthesis)  for $h=2^{-6}$ and different values of $\nu$ and $\omega$.\label{Tbl2}}
	\begin{center}
		\label{exact}
		\scalebox{.8}
		{
\begin{tabular}{|c|c|c|c|c|c|c|c|c|c|c|c} \hline
&$\nu\setminus\omega$ &$10^{-4}$ &$10^{-3}$&$10^{-2}$ & $10^{-1}$ & $1$ &$10^{1}$&$10^{2}$&$10^{3}$&$10^{4}$\\ \hline

				IASSS & $10^{-2}$& 56(0.59) & 56(0.55) & 56(0.55) & 56(0.57) & 56(0.57) & 55(0.54) & 50(0.42)  & 40(0.31) & 51(0.29) \\
				& $10^{-4}$& 50(0.45) & 50(0.42) & 50(0.41) & 50(0.41) & 50(0.41) & 50(0.41) & 48(0.39)  & 40(0.32) & 51(0.30) \\ 
				& $10^{-6}$& 40(0.34) & 40(0.34) & 40(0.31) & 40(0.31) & 40(0.30) & 40(0.31) & 40(0.31)  & 42(0.30) & 51(0.29) \\
				& $10^{-8}$& 51(0.32) & 51(0.29) & 51(0.30) & 51(0.29) & 51(0.30) & 51(0.30) & 51(0.29)  & 51(0.31) & 52(0.30) \\
				\hline

IBAS  & $10^{-2}$& 38(0.36) & 38(0.34) & 38(0.33) & 38(0.32) & 38(0.33) & 24(0.23) & 467(2.52) & $\dag$ & $\dag$ \\
	 & $10^{-4}$& 35(0.27) & 35(0.25) & 35(0.25) & 33(0.25) & 35(0.25) & 35(0.25) & 39(0.28)  & $\dag$ & $\dag$ \\ 
	 & $10^{-6}$& 33(0.23) & 33(0.21) & 33(0.20) & 33(0.20) & 33(0.21) & 33(0.20) & 33(0.20)  & 56(0.31) & $\dag$ \\
	 & $10^{-8}$& 38(0.23) & 38(0.21) & 38(0.22) & 38(0.21) & 38(0.21) & 38(0.20) & 38(0.20)  & 38(0.20) & 65(0.28) \\
				
				\hline\hline\hline
				
P-ASSS& $10^{-2}$& 28(0.44) & 28(0.41) & 28(0.40) & 28(0.38) & 28(0.38) & 29(0.42) & 29(0.35) & 21(0.25) & 21(0.16) \\
       & $10^{-4}$& 29(0.41) & 29(0.39) & 29(0.38) & 29(0.39) & 29(0.38) & 29(0.38) & 29(0.36) & 25(0.24) & 21(0.16)\\ 
	   & $10^{-6}$& 25(0.30) & 25(0.26) & 25(0.27) & 25(0.26) & 25(0.25) & 25(0.24) & 25(0.25) & 25(0.23) & 21(0.16) \\
	   & $10^{-8}$& 21(0.25) & 21(0.18) & 21(0.17) & 21(0.17) & 21(0.20) & 21(0.17) & 21(0.18) & 21(0.16) & 21(0.16) \\
		\hline
P-BAS & $10^{-2}$& 18(0.27) & 19(0.31) & 20(0.32) & 20(0.32) & 20(0.32) & 17(0.25) & 26(0.30) & 49(0.45) & 44(0.35) \\
& $10^{-4}$& 20(0.26) & 20(0.25) & 21(0.26) & 22(0.29) & 22(0.25) & 22(0.27) & 20(0.20) & 47(0.44) & 48(0.35)\\ 
& $10^{-6}$& 18(0.17) & 19(0.19) & 20(0.19) & 21(0.23) & 21(0.21) & 21(0.19) & 22(0.21) & 28(0.21) & 47(0.37) \\
& $10^{-8}$& 17(0.13) & 19(0.16) & 20(0.16) & 20(0.17) & 21(0.17) & 21(0.19) & 22(0.18) & 22(0.20) & 28(0.20) \\
\hline

P-PRESB & $10^{-2}$& 7(0.31) & 7(0.20) & 7(0.29) & 8(0.40) & 10(0.62) & 23(2.92) & 107(33.19) & $\ddag$   & $\ddag$ \\
        & $10^{-4}$& 8(0.28) & 8(0.17) & 8(0.20) & 8(0.26) &  9(0.32) & 12(0.76) & 60(10.00) & $\ddag$   & $\ddag$\\ 
        & $10^{-6}$& 8(0.24) & 8(0.14) & 8(0.12) & 8(0.19) &  8(0.21) &  9(0.24) & 12(0.55) & 91(14.59) & $\ddag$ \\
        & $10^{-8}$& 8(0.22) & 8(0.10) & 8(0.11) & 8(0.11) &  8(0.17) &  8(0.18) & 9(0.17) & 12(0.14)  & 81(11.82) \\
\hline
P-BD & $10^{-2}$& 14(0.33) & 14(0.25) & 14(0.21) & 16(0.23) & 22(0.36) & 54(0.81) & 182(4.42) & 174(3.86) & 84(1.16) \\
     & $10^{-4}$& 16(0.30) & 16(0.21) & 16(0.20) & 22(0.26) & 42(0.54) &195(5.09) & $\ddag$   & $\ddag$   & 154(3.02)\\ 
     & $10^{-6}$& 16(0.23) & 16(0.15) & 16(0.14) & 17(0.16) & 38(0.35) &185(4.25) & $\ddag$   & $\ddag$   & 168(3.31)\\
     & $10^{-8}$& 15(0.18) & 15(0.12) & 15(0.10) & 15(0.09) & 22(0.17) & 85(1.07) & 153(2.97) & 166(3.27) & 116(1.80) \\
\hline

\end{tabular}}
\end{center}
\end{table}

\begin{table}[!h]
	\centering\caption{Number of iterations of the methods along with the elapsed CPU time (in parenthesis)  for $h=2^{-7}$ and different values of $\nu$ and $\omega$.\label{Tbl3}}
	\begin{center}
		\label{exact}
		\scalebox{.8}
		{
\begin{tabular}{|c|c|c|c|c|c|c|c|c|c|c|c} \hline
&$\nu\setminus\omega$ &$10^{-4}$ &$10^{-3}$&$10^{-2}$ & $10^{-1}$ & $1$ &$10^{1}$&$10^{2}$&$10^{3}$&$10^{4}$\\ \hline
			IASSS& $10^{-2}$& 57(3.43) & 57(3.30) & 57(3.26) & 57(3.35) & 57(3.31) & 56(2.95) & 52(2.45)  & 44(1.70) & 51(1.41) \\
			& $10^{-4}$& 52(2.46) & 52(2.36) & 52(2.43) & 52(2.39) & 52(2.43) & 52(2.56) & 52(2.48)  & 44(1.78) & 51(1.42) \\ 
			& $10^{-6}$& 44(1.78) & 44(1.65) & 44(1.67) & 44(1.78) & 44(1.78) & 44(1.79) & 44(1.73)  & 43(1.60) & 51(1.44) \\
			& $10^{-8}$& 51(1.53) & 51(1.49) & 51(1.51) & 51(1.42) & 51(1.43) & 51(1.45) & 51(1.50)  & 51(1.45) & 52(1.35) \\
			\hline

IBAS & $10^{-2}$& $\dag$   & $\dag$   & $\dag$   & $\dag$   & 48(2.30) & 25(1.32) & 465(14.04) & $\dag$ & $\dag$ \\
     & $10^{-4}$& 36(1.43) & 36(1.43) & 36(1.43) & 36(1.46) & 36(1.40) & 36(1.42) & 39(1.70)   & $\dag$ & $\dag$ \\ 
	 & $10^{-6}$& 33(1.09) & 33(1.10) & 33(1.08) & 33(1.10) & 33(1.10) & 33(1.14) & 33(1.13)   & 56(1.78) & $\dag$ \\
	 & $10^{-8}$& 38(0.98) & 38(0.98) & 38(0.97) & 38(0.97) & 38(0.96) & 38(0.96) & 38(0.97)   & 38(0.97) & 64(1.48) \\
				\hline\hline\hline

P-ASSS& $10^{-2}$& 28(2.40) & 28(2.40) & 28(2.40) & 28(2.40) & 28(2.40) & 29(2.42) & 30(2.42) & 28(1.40) & 25(0.92) \\
      & $10^{-4}$& 30(2.38) & 30(2.42) & 30(2.44) & 30(2.38) & 30(2.54) & 30(2.37) & 29(2.23) & 28(1.51) & 25(0.98)\\ 
	  & $10^{-6}$& 28(1.55) & 28(1.59) & 28(1.57) & 28(1.52) & 28(1.50) & 28(1.49) & 28(1.48) & 27(1.26) & 25(0.95) \\
	  & $10^{-8}$& 25(0.95) & 25(0.94) & 25(0.94) & 25(0.94) & 25(0.94) & 25(0.94) & 25(0.92) & 25(0.90) & 24(0.79) \\
	\hline
P-BAS & $10^{-2}$& 18(1.77) & 19(1.84) & 20(1.87) & 20(1.88) & 20(1.83) & 17(1.61) & 26(1.77) & 49(2.33) & 46(1.47)\\
& $10^{-4}$& 20(1.35) & 21(1.42) & 22(1.41) & 22(1.50) & 22(1.47) & 22(1.54) & 20(1.30) & 47(2.15) & 50(1.71)\\ 
& $10^{-6}$& 18(0.76) & 19(0.77) & 20(0.80) & 21(0.90) & 21(0.81) & 21(0.87) & 22(0.96) & 28(1.16) & 49(1.67) \\
& $10^{-8}$& 18(0.48) & 19(0.52) & 20(0.56) & 30(0.53) & 21(0.59) & 22(0.62) & 22(0.64) & 22(0.59) & 28(0.81) \\
\hline

P-PRESB & $10^{-2}$& 7(1.15) & 7(1.09) & 7(1.92) & 8(2.48) & 10(3.79) & 23(18.42) & $\ddag$   & $\ddag$   & $\ddag$ \\
        & $10^{-4}$& 8(1.04) & 8(0.92) & 8(1.25) & 8(1.62) &  9(1.73) & 12(4.16)  & 60(58.46) & $\ddag$   & $\ddag$\\ 
        & $10^{-6}$& 8(0.68) & 8(0.60) & 8(0.58) & 8(0.95) &  8(0.97) &  9(1.34)  & 12(2.56)  & 92(70.74) & $\ddag$ \\
        & $10^{-8}$& 8(0.51) & 8(0.41) & 8(0.40) & 8(0.40) &  8(0.65) &  8(0.62)  & 9(0.81)   & 12(1.66)  & 87(47.56) \\

\hline

P-BD & $10^{-2}$& 14(1.45) & 14(1.33) & 14(1.35) & 16(1.50) & 22(1.90) & 54(4.12)  & 182(17.19) & 188(15.20) & 142(9.10) \\
     & $10^{-4}$& 16(1.21) & 16(1.24) & 18(1.27) & 22(1.51) & 42(2.76) &199(19.97) & $\ddag$    & $\ddag$    & $\ddag$\\ 
     & $10^{-6}$& 16(0.75) & 16(0.67) & 16(0.66) & 20(0.79) & 40(1.64) &203(17.87) & $\ddag$    & $\ddag$    & $\ddag$\\
     & $10^{-8}$& 15(0.47) & 15(0.41) & 15(0.40) & 16(0.42) & 30(0.89) &151(10.12) & $\ddag$    & $\ddag$    & $\dag$ \\
\hline

			\end{tabular}}
		\end{center}
	\end{table}

\section{Conclusion}\label{Sec6}

 We have presented a method called  alternating symmetric positive definite and  scaled symmetric positive semidefinite splitting (ASSS) method, for solving the system arisen from finite element discretization of a distributed optimal control problem with time-periodic parabolic equations. We have proved that the method is unconditionally  convergent. We have compared the numerical results of the ASSS method and the corresponding induced preconditioner with those of the BAS iteration method. We also compared the numerical results of the ASSS preconditioner to those of the PRESB preconditioner. Numerical results showed that the proposed \cred{method has some advantages over the two other tested methods. }

\section*{Acknowledgments}

The author would like to thank the referees for their careful reading of the paper and giving several valuable comments and suggestions.


\end{document}